\documentclass[12pt]{article}
\usepackage{latexsym}
\usepackage{amsthm}
\usepackage{amsmath}
\usepackage{epsfig}
\usepackage{amsfonts}
\usepackage{amssymb}
\newtheorem{theorem}{Theorem}
\newtheorem{lemma}{Lemma}[theorem]

\usepackage{graphics}

\def\R{\mathbb R}

\begin{document}

\title{Knots of Constant Curvature}

\author{Jenelle Marie McAtee \footnote{Department of Mathematics, University of Iowa, Iowa
City IA 52240.  This work is part of the author's ongoing PhD thesis research
conducted under the guidance of thesis advisor Jonathan Simon.  The author has
been partially supported by NSF grants DMS9706789 and DMS0107209. }}

\maketitle

\markboth{J. McAtee}{Knots of Constant Curvature}

\begin{abstract} In this paper we show how to realize all knot (and link) types
as $\mathcal{C}^{2}$ smooth curves of constant curvature.  Our proof is
constructive: we build the knots with copies of a fixed finite number of
"building blocks" that are particular segments of helices and circles.  We use
these building blocks to construct all closed braids.
\end{abstract}

\section{Introduction}  Circles and helices are standard examples of smooth
curves of constant curvature.  To construct other $\mathcal{C}^{2}$ curves with
constant curvature we may splice together pieces of helices and circles with the
same curvature in such a way that the resulting curve is $\mathcal{C}^{2}$ (see
below).  In \cite{EK}, the authors Koch and Engelhardt integrated piecewise circular curves on
$S^{2}$ to obtain nonplanar unknots of constant curvature.  By manipulating
several helix and circle segments, we found the granny knot shown in figure one.  In this paper we
develop a systematic method for constructing knots and show that all knot and
link types can be realized as $\mathcal{C}^{2}$ smooth curves with constant
curvature.  In a subsequent paper, we develop a different technique for realizing
many $\mathcal{C}^{2}$ constant curvature knots using the method of Koch and
Engelhardt.

\section{The splicing method}

Suppose that $\Gamma_{1}$ and $\Gamma_{2}$ are $\mathcal{C}^{2}$
smooth curves of the same constant
nonvanishing curvature $\kappa$. Let $r_{1}: [0,\mathcal{L}_{1}] \rightarrow $
$\R^{3}$ and $r_{2}: [0,\mathcal{L}_{2}] \rightarrow \R^{3}$ be $\mathcal{C}^{2}$ arclength parameterizations of
$\Gamma_{1}$ and $\Gamma_{2}$ respectively.  For $i \in \{1,2\}$, let
$\{T_{i}(*),N_{i}(*),B_{i}(*)\}$ denote the Frenet frame of $r_{i}$ at the point
$r_{i}(*).$ 
Since the Frenet frame is an orthonormal set of vectors, there exists a rotation $A:\R^{3}
\rightarrow \R^{3}$ such that $A(T_{2}(0)) = T_{1}(\mathcal{L}_{1}), A(N_{2}(0)) =
N_{1}(\mathcal{L}_{1}),$ and $A(B_{2}(0)) = B_{1}(\mathcal{L}_{1}).$We will define
$u(s)$ to be the curve obtained by rotating $\Gamma_{2}$ so that its initial
Frenet frame coincides with the terminal Frenet frame of $\Gamma_{1}$ and
translating the rotated curve so that its initial point is $r_{1}(L_{1})$. In
particular let $w =
r_{1}(\mathcal{L}_{1}) - A(r_{2}(0))$ and define $u(s)$ by $u(s) =
A(r_{2}(s-\mathcal{L}_{1})) + w$
for $s \in [\mathcal{L}_{1}, \mathcal{L}_{1}+\mathcal{L}_{2}]$.
 Consider the curve $\Gamma_{1}*\Gamma_{2}$ given by
the parameterization \begin{equation}
                             r(s) = \nonumber
			       \begin{cases}
			          r_{1}(s)  & \text{if $s\in [0,\mathcal{L}_{1}]$}, \\
				  u(s) & \text{if $s\in
				  [\mathcal{L}_{1}, \mathcal{L}_{1}+\mathcal{L}_{2}]$}
				\end{cases}
		      \end{equation}
where $s$ is the arclength parameter of $\Gamma_{1}$.  For a discussion of Frenet
frames and curvature please see \cite{DoC}.

Since $A$ is linear, it follows that
$\frac{du}{ds}(s) = A(r'_{2}(s-\mathcal{L}_{1}))$, and consequently $u'(\mathcal{L}_{1}) =
A(r'_{2}(\mathcal{L}_{1}-\mathcal{L}_{1})) = A(r'_{2}(0)) = A(T_{2}(0)) = T(\mathcal{L}_{1}) =
r'_{1}(\mathcal{L}_{1})$.  Since $r$ is defined piecewise by the $\mathcal{C}^{2}$ functions
$r_{1} \text{and } u$, it follows from this that $r$ is differentiable and $\mathcal{C}^{1}$.
Similarly since $\frac{d^2 u}{ds^2}(s) = A(r''_{2}(s-\mathcal{L}_{1})) =
A(T'_{2}(s-\mathcal{L}_{1})) = A(\kappa N_{2}(s-\mathcal{L}_{1})) =
\kappa A(N_{2}(s-\mathcal{L}_{1}))$ and since $u''(\mathcal{L}_{1}) =
\kappa A(N_{2}(\mathcal{L}_{1}-\mathcal{L}_{1})) = \kappa A(N_{2}(0)) =
\kappa N_{1}(\mathcal{L}_{1}) = T'_{1}(\mathcal{L}_{1}) = r''_{1}(\mathcal{L}_{1})$, it follows
that $r$ is $\mathcal{C}^{2}$. Furthermore since $r$ is an arclength parameterization of
$\Gamma_{1}*\Gamma_{2}$ and $r''$ is given by \begin{equation}
                             r''(s) = \nonumber
			       \begin{cases}
			          \kappa N_{1}(s)  & \text{if $s\in [0,\mathcal{L}_{1}]$}, \\
				  \kappa A(N_{2}(s-\mathcal{L}_{1})) & \text{if $s\in
				  [\mathcal{L}_{1}, \mathcal{L}_{1}+\mathcal{L}_{2}]$}
				\end{cases}
		      \end{equation}
it follows that $\Gamma_{1}*\Gamma_{2}$ has constant curvature $\kappa$.

Essentially $\Gamma_{1}*\Gamma_{2}$ is obtained by rotating $\Gamma_{2}$ so that its initial
Frenet frame coincides with the terminal Frenet frame of $\Gamma_{1}$ and then attaching the
beginning of the rotated $\Gamma_{2}$ to the end of $\Gamma_{1}$.  This process of joining $\Gamma_{2}$
to $\Gamma_{1}$ to achieve $\Gamma_{1}*\Gamma_{2}$, a $\mathcal{C}^{2}$ smooth curve, will be refered to as
\emph{splicing} $\Gamma_{2}$ to $\Gamma_{1}$.  Note that the splicing operation is
associative. Hence for the sake of simplicity, $\Gamma_{1}\Gamma_{2}$ will denote
$\Gamma_{1}*\Gamma_{2}$ and parenthesis will be ignored.

The curvature $\kappa$ of a helix parameterized by $\alpha (t) = [\pm rcos(t),\pm
rsin(t),\pm ht] $ is given by $\kappa = \frac{r}{r^{2}+h^{2}}$.  Hence for any fixed
$\kappa \neq 0 $, there is an infinite number of parameter values for $r$ and $h$
which determine a helix with constant curvature $\kappa$.  Consequently there are
infinitely many helices which can be used to form $\mathcal{C}^{2}$ curves of
constant curvature $\kappa$ (see figure 1). In this paper, we focus on two curves
of curvature one: the unit circle and the helix with $r=h=\frac{1}{2}.$
\begin{figure}[here]
\centering 
\includegraphics[height=4in]{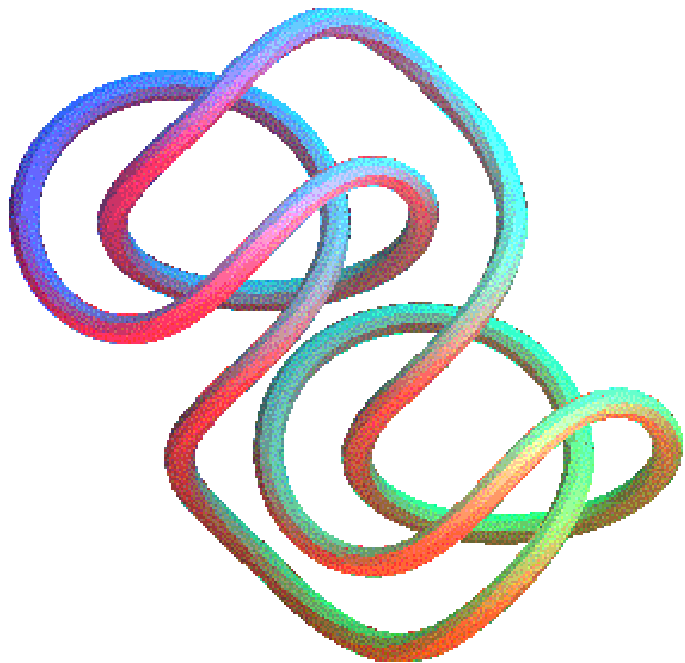}
\caption{a $\mathcal{C}^{2}$ granny knot with constant curvature constructed
by splicing together pieces of helicies}
\label{Figure 1:}
\end{figure}

\section{Fundamental Building Blocks}

Let the vector notation $[r,h,t_{i},t_{e}]$ denote the piece of the helix  $\alpha
(t) = [rcos(t),rsin(t),ht] $ with initial point $\alpha (t_{i})$ and terminal point $\alpha
(t_{e})$.  We define the following \emph{elementary pieces}: \begin{eqnarray}
\nonumber a & = &
[\frac{1}{2}, \frac{1}{2}, 0, \frac{\pi}{2}] \\ \nonumber b & = & [1, 0, 0, \pi]
\\ \nonumber c & = & [\frac{1}{2},\frac{1}{2}, 0, 2\pi] \\ \nonumber d & = & [\frac{1}{2}, \frac{-1}{2},
\frac{3\pi}{2}, 2\pi] \\ \nonumber e & = &  [\frac{1}{2}, \frac{-1}{2}, 0, 2\pi]
\\ \nonumber f & = & [\frac{1}{2}, \frac{-1}{2}, 0, \frac{3\pi}{2}] \\ \nonumber g & = & [\frac{1}{2}, \frac{-1}{2}, 0,
\frac{\pi}{2}] \\ \nonumber l & = & [\frac{1}{2},\frac{1}{2}, 0, 8\pi] \end{eqnarray} 
All of these helices have curvature equal to one.  In fact, the elementary pieces
are either part of the unit circle or part of a helix with $r = h = \frac{1}{2}$. 
 After reparametrizing by arclength we may splice together the elementary pieces to
 form the following family of curves 
 $ \mathcal{B} =
\{i_{+},i_{-},j_{+},j_{-},k_{+},k_{-}\}$ where $i_{+} = (abcd)^{4}$, $i_{-} =
(aebd)^4$, $j_{+} = (adbe)^4$, $j_{-} = (adbc)^4$, $k_{+} = l$, and $k_{-} =
abfdgbabfdfgb$.  We will call elements of $\mathcal{B}$ \emph{sticks}.  

Let $F$ denote the Frenet frame at the initial point of
$a$.  The sticks are constructed so that, as parameterized curves, the Frenet frame at
the initial and terminal point of each stick is $F$.  This property affects the
geometry of the sticks.  For example, if we subtract the initial point of $k_{-}$
from its terminal point, we get the displacement vector $[0,0,-4\pi]$.  Hence $k_{-}$ accomplishes a net
movement in the negative $z$-direction.  However since we require the Frenet frame of
$k_{-}$ to be $F$ at its initial and terminal points and since the $z$-coordinate of the
unit tangent in $F$ is positive, $k_{-}$ must travel in the positive $z$-direction
before turning downwards (see figure 2). Similarly, each of the remaining sticks
must twist about in space in order for the initial and terminal Frenet frame to be
$F$.  Because of this twisting and because we will eventually construct knots, we
have to be careful to make sure that the spliced sticks are simple curves (see
figures 2 and 3).  Furthermore, since every element of $\mathcal{B}$ has its initial and terminal Frenet frame
equal to $F$, any two sticks can be spliced together without rotation.  Hence the
resulting curve also has the attribute that its initial and terminal Frenet frame
is $F$.
\pagebreak
\begin{figure}
\centering 
\includegraphics[height=3in]{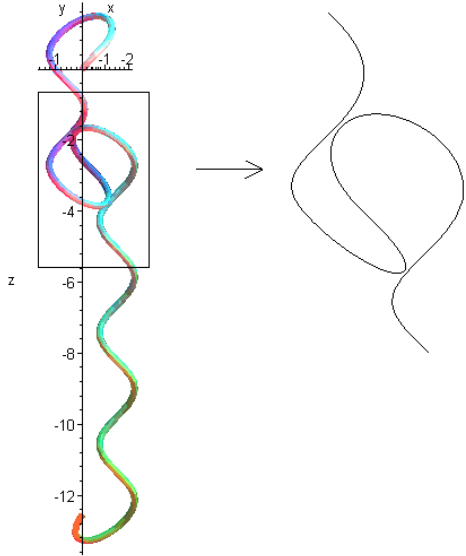}
\caption{The stick $k_{-}$ has been thickened so it can be viewed more easily. 
The piece of the curve in the box has been redrawn as a spacecurve and magnified
to demonstrate that the stick $k_{-}$ does not intersect itself.}
\label{Figure 2: $k_{-}$}
\end{figure}

In addition to requiring that every stick have $F$ as its initial and terminal
Frenet frame, the sticks have been constructed so that the distance from the
initial point of each stick to its terminal point is the same. If we call the vector from the initial point of the stick to its
terminal point the \emph{displacement vector} of the stick, then the displacement
vectors of $i_{+}$, $i_{-}$, $j_{+}$, $j_{-}$, $k_{+}$, and $k_{-}$ are
$(4\pi,0,0)$, $(-4\pi,0,0)$, $(0,4\pi,0)$, $(0,-4\pi,0)$, $(0,0,4\pi)$, and
$(0,0,-4\pi)$ respectively.  Therefore since sticks and curves created from
elements of $\mathcal{B}$ are spliced by simply translating the
beginning of one stick to the end of another, it follows that given any path on a
square lattice
where the unit length is $4\pi$, we can imitate the path with a
$\mathcal{C}^{2}$ curve of constant curvature by splicing together elements of
$\mathcal{B}$.  

\section{Knots}
In this section we will show how to splice together elementary pieces to realize a
knot type as a $\mathcal{C}^{2}$ curve of constant curvature.  We will realize the
knot types by forming braid closures.

While any number of elementary pieces can be spliced together to form a
$\mathcal{C}^{2}$ curve of constant curvature, not every such curve is simple. 
Each stick has been constructed to be a simple curve.  Figure 3 shows the sticks
$i_{+}$, $i_{-}$, $j_{+}$, and $j_{-}$.  A view of $k_{-}$ is shown in figure 2
along with enlargements of parts of the curve to demonstrate that there are indeed
no self-intersections.  The stick $k_{+}$ is simply a piece of a single helix, and
therefore clearly does not intersect itself.
\begin{figure}[here]
\centering
\includegraphics[height=4in,clip=]{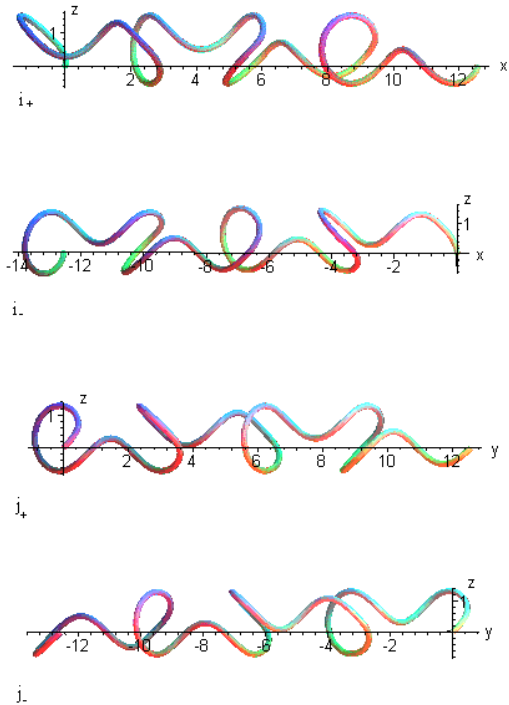}
\caption{The sticks have all been thickened so they can be viewed more easily.  If
the sticks are rotated in $\R^{3}$, then it is clear that they are simple curves.}
\label{Figure 3: Sticks}
\end{figure}

Though each stick in $\mathcal{B}$ is simple, not every word made from elements of
$\mathcal{B}$ represents a simple curve.  For example, $i_{+}i_{-}$, $i_{-}i_{+},$
$j_{+}j_{-},$ $j_{-}j_{+},$ $k_{+}k_{-},$ $k_{-}k_{+},$ $k_{-}j_{+},$ and
$i_{-}k_{-}$ are all curves with self intersections.  We will consider these
letter pairs as \emph{not allowable} because they lead to self-intersection.  All
other letter pairs represent simple curves and consequently these remaining letter
pairs will be called \emph{allowable letter pairs}.  As an example of a curve
represented by an allowable letter pair, $k_{-}k_{-}$ is shown in figure 4.     

\begin{figure}[here]
\centering
\includegraphics[height=3in,clip=]{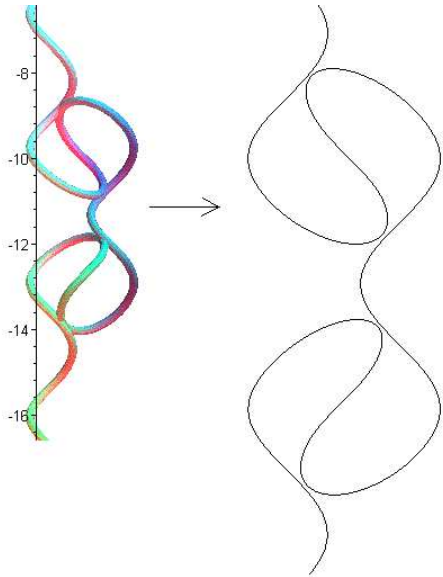}
\caption{ Here we have close views of the join after a $k_{-}$ stick has been
spliced to another $k_{-}$ stick. In the first drawing, the curve has been
thickened for ease of viewing and this thickening accounts for the
self-intersections.  The second drawing shows the join as a spacecurve, and it is
clear that after the join the curve remains simple.}
\label{Figure 4: $k_{-}k_{-}$}
\end{figure}
Let $I_{+}$, $I_{-}$, $J_{+}$, $J_{-}$, $K_{+}$, and $K_{-}$ denote $(4\pi,0,0)$, $(-4\pi,0,0)$, $(0,4\pi,0)$, $(0,-4\pi,0)$, $(0,0,4\pi)$, and
$(0,0,-4\pi)$ respectively.  We call any word in $I_{\pm}$, $J_{\pm}$, and $K_{\pm}$  a
\emph{lattice polygon}, and we identify the word with an actual polygon on a lattice
(with unit length equal to $4\pi$) in $\R^{3}$.  As before we will call  $I_{+}I_{-}$, $I_{-}I_{+},$
$J_{+}J_{-},$ $J_{-}J_{+},$ $K_{+}K_{-},$ $K_{-}K_{+},$ $K_{-}J_{+},$ and
$I_{-}K_{-}$ \emph{unallowable letter pairs}.  We have the following lemma:

\begin{lemma} If $P$ is a lattice polygon without self-intersection and without
unallowable letter pairs, then the same word written with elements of $\mathcal{B}$
represents a simple curve.
\end{lemma}

\begin{proof}  Let $P = E_{1}E_{2}...E_{n}$ where
$E_{m}\in\{I_{\pm},J_{\pm},K_{\pm}\}$, and let $p$ = $e_{1}e_{2}...e_{n}$ where 
$e_{m}\in\{i_{\pm},j_{\pm},j_{\pm}\}$, $e_{m}$ is the lowercase of $E_{m}$, and
$e_{m}$ has the same subscript as $E_{m}$.  Then $p$ realizes the polygon $P$ as a
$\mathcal{C}^{2}$ curve of constant curvature.  Since $P$ does not contain any
unallowable letter pairs, $p$ also does not contain any unallowable letter pairs. 
Therefore adjacent sticks in $p$ do not intersect.  It is easy to show that each
stick $i_{\pm}$, $j_{\pm}$, $k_{\pm}$ is contained in a tube of radius $\pi$
centered about the line containing the initial and terminal points of the stick. 
Therefore since $P$ is simple and since nonconsecutive sticks in $p$ are separated
by a distance of at least $4\pi$, it follows that the nonconsecutive sticks in $p$
do not intersect.  Thus we have that $p$ is a simple curve.
\end{proof}

\subsection{Braids} 

Braid components can be formed from the sticks in such a way that if each strand
in the braid component is expressed as a word in elements of $\mathcal{B}$, the
words contain only allowable word pairs.  It follows from the previous lemma that
each strand in the braid component is a simple $\mathcal{C}^{2}$ curve of constant curvature.  In
figure 5, we have views of a braid component that contains a positive
crossing.  The $n^{th}$ strand is formed by translating the curve
$k_{+}k_{+}i_{+}k_{+}k_{+}k_{+}$ so that its initial point is $(-4\pi,0,0)$.  We translate the
curve without rotation so that the Frenet frame at the initial and terminal points
of the translated curve is still $F$. The $(n+1)^{th}$ strand is formed by
translating (without rotation) the curve
$k_{+}j_{-}k_{+}k_{+}i_{-}k_{+}j_{-}k_{+}$ so that its initial point is the origin.  The remaining strands
are simply the curve $(k_{+})^{5}$ translated (without rotation) so that the initial
point of each strand is an integer
multiple of $4\pi$. Views of a braid component containing a negative crossing are given in figure 6. 
The $n^{th}$ strand is the curve $k_{+}j_{+}k_{+}i_{+}k_{+}j_{-}k_{+}$ translated
(without rotation) so that its initial point is $(4\pi,0,0)$, and the $(n+1)^{th}$ strand is the
curve $k_{+}k_{+}i_{-}k_{+}k_{+}$ translated (without rotation) to have its initial
point at the origin.  As
before the remaining strands are the curves $(k_{+})^{5}$ translated so that the
inital points are all integer
multiples of $4\pi$.
\pagebreak
\begin{figure}[here]
\centering
\includegraphics[height=4in,clip=]{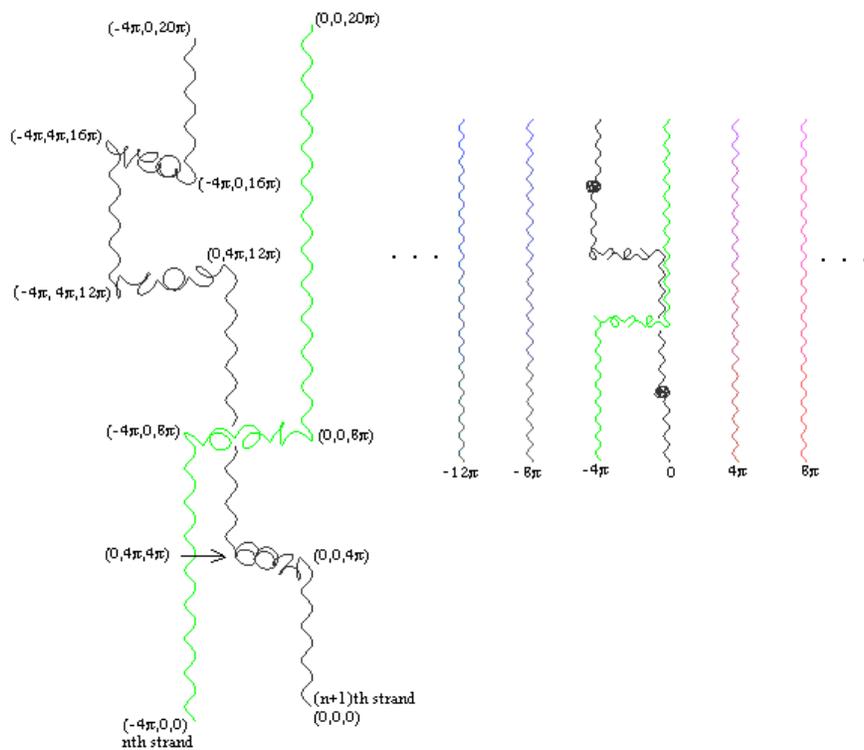}
\caption{a braid component with a positive crossing}
\label{Figure 5: }
\end{figure}

\begin{figure}[here]
\centering
\includegraphics[height=4in,clip=]{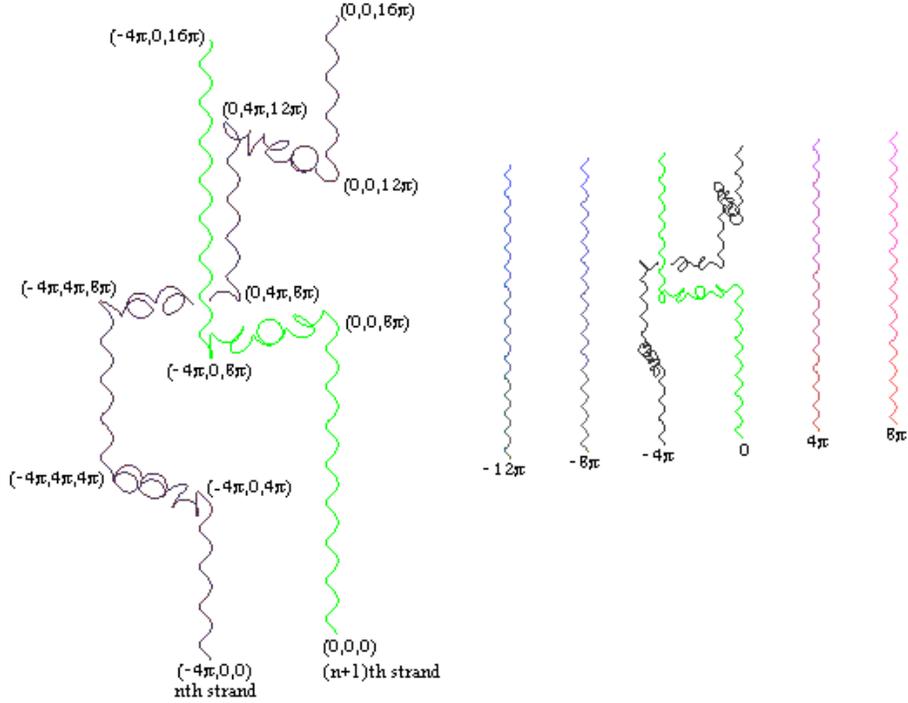}
\caption{a braid component with a negative crossing}  
\label{Figure 6: }
\end{figure}
\pagebreak
The braid components have been constructed so that the initial and terminal points of
each strand are on equally spaced lines contained in the $(x,z)$-plane and
perpendicular to the $x$-axis.  This allows the braid components to be nicely
stacked.  Furthermore since the braid components are constructed from elements of
$\mathcal{B}$, the Frenet frame at the initial and terminal points of each strand is $F$. 
Hence the braid components can be stacked in a $\mathcal{C}^{2}$ fit, and each strand in
the resulting braid remains a $\mathcal{C}^{2}$ curve of constant curvature (see
figure 7).
\begin{figure}[here]
\centering
\includegraphics[height=4in]{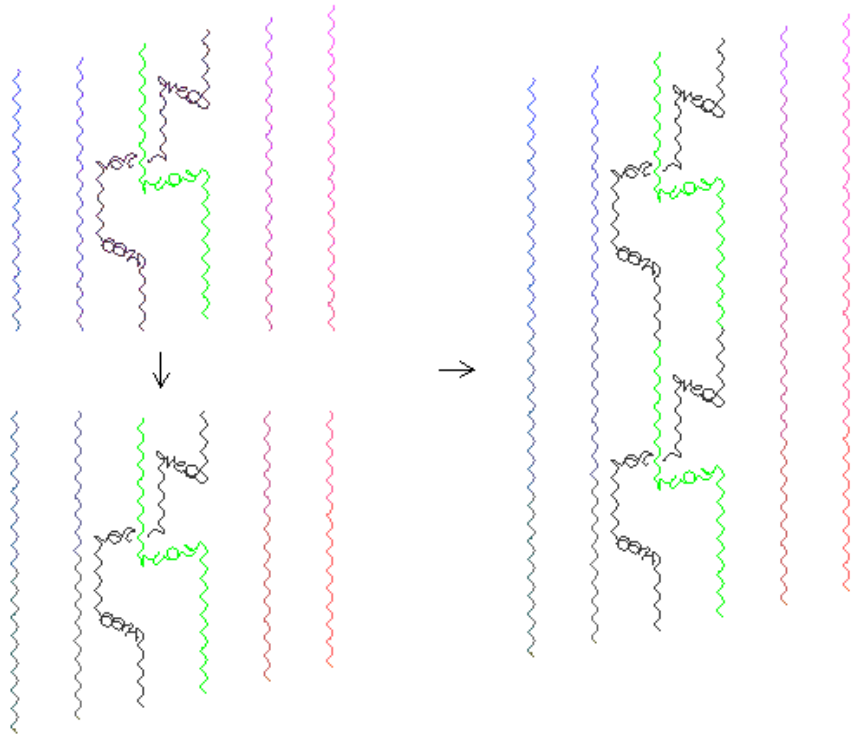}
\caption{A six strand braid component containing a negative
crossing is spliced to the top of another braid component containing a negative
crossing to form a $\mathcal{C}^{2}$ braid with two crossings.}
\label{Figure 7: }
\end{figure}
\pagebreak
\subsection{Closing pieces}

In addition to the braid components, \emph{closing pieces} may also be fashioned from
elements of $\mathcal{B}$.  Again these closing pieces are constructed so that when they
are expressed as words made from the elements of $\mathcal{B}$, only allowable letter pairs
appear.  In figure 8 we have constructed a closing piece.  This is the curve
$k_{+}i_{+}k_{-}k_{-}k_{-}k_{-}k_{-}i_{+}k_{+}$ which has been translated without
rotation so that its inital point is $(0,0,12\pi)$.  By adding and deleting an
appropriate number of sticks, the lengths in the various directions of the closing
pieces can be adjusted.
\pagebreak   
\begin{figure}[here]
\centering
\includegraphics[height=5in]{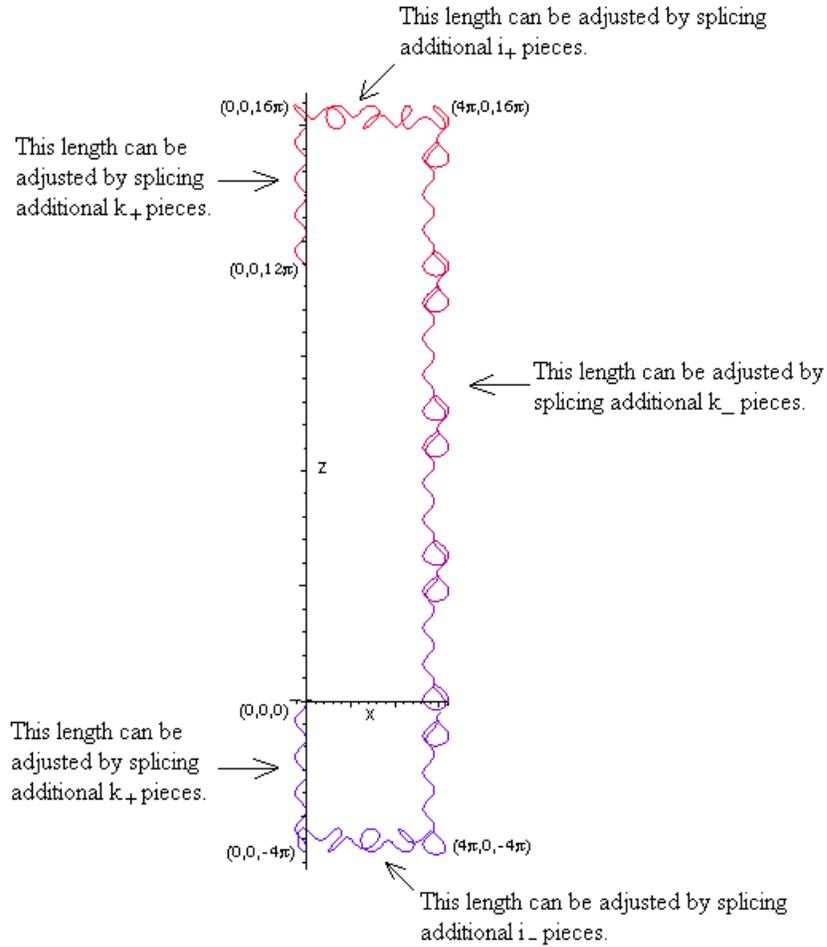}
\caption{an adjustable closing piece}
\label{Figure 8: }
\end{figure}

Since we can adjust the dimensions of the closing pices and since the initial and terminal
points of each strand are on equally spaced lines contained in the $(x,z)$-plane and
perpendicular to the $x$-axis, the closing pieces can be attached to the braids to form the
braid closure.  As before since the closing pieces are formed from elements of
$\mathcal{B}$, the initial and terminal Frenet frames of the closure pieces are $F$, and
consequently the closure of the braid is a $\mathcal{C}^{2}$ curve of constant curvature. 
Since Alexander's Theorem (see \cite{Mu}) states that every knot (and link) can be realized as the closure of a braid, we
have the following theorem: 
\pagebreak
\begin{theorem}  Every knot (and link) can be represented by a $\mathcal{C}^{2}$ curve of constant curvature.
\end{theorem} 

In figure 9 we show a $\mathcal{C}^{2}$ trefoil of constant curvature constructed
from sticks as the closure of a braid.
\begin{figure}[here]
\centering 
\includegraphics[height=5in,width = 5in]{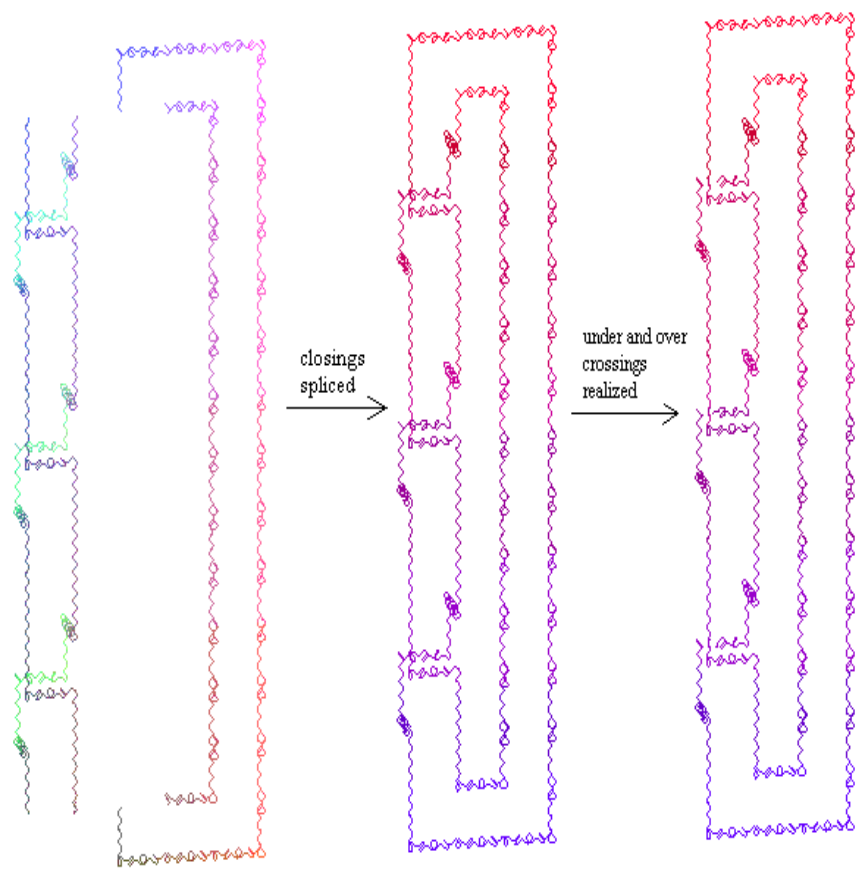}
\caption{a $\mathcal{C}^{2}$ trefoil of constant curvature}
\label{Figure 9: }
\end{figure}

\newpage

\bibliographystyle{siam}

\bibliography{Knot_Constant_Curvature}

\begin{thebibliography}{1}

\bibitem{DoC}
{\sc M.~D. Carmo}, {\em Differential Geometry of Curves and Surfaces},
  Prentice-Hall, New Jersey, 1976.

\bibitem{EK}
{\sc C.~Engelhardt and R.~Koch}, {\em Closed space curves of constant curvature
  consisting of arcs of circular helices}, Journal for Geometry and Graphics, 2
  (1998), pp.~17--31.

\bibitem{Mu}
{\sc K.~Murasugi}, {\em Knot Theory and Its Applications}, Burkhauser, Boston,
  1996.

\end{thebibliography}

\end{document}